\input amstex
\input amsppt1
\loadmsbm
\UseAMSsymbols
\NoBlackBoxes
\parindent=0pt
\magnification=1200
\def\cd{\Cal{D}}
\def\tc{\witi C}

 sCaled 850

\def\td{\widetilde{\Cal D}}
\def\cd{\Cal D}

\def\whh{\widehat l}

\def\enpr{\quad \vrule height .9ex width .8ex depth -.1ex}

\def\bar{\overline}

\def\bre{\Bbb R}

\def\llim{\mathop{\longrightarrow}}

\def\vo{V^{(0)}}
\def\von{V^{(1)}}

\def\rvo{\bre^{V^{(0)}}}
\def\rvon{\bre^{V^{(1)}}}

\def\Lb{\Lambda}

\def\cl{{\Cal L}}
\def\cf{{\Cal F}}

\def\s,{\quad $\,$}

\def\qua{\quad$\,$}
\def\smad{\smallskip Proof.\ }

\def\s,{\quad $\,$}

\def\witi{\widetilde}

\def\disp{\displaystyle}

\bigskip

 \centerline{\bf A P.C.F. Self-Similar Set with no Self-Similar Energy}
 
\bigskip
\centerline{\bf Roberto Peirone}
\bigskip
\centerline{Universit\`a di Roma "Tor Vergata", Dipartimento di
Matematica}

\centerline{via della Ricerca Scientifica, 00133, Roma, Italy}
\smallskip
\centerline{e.mail: peirone\@mat.uniroma2.it}

\bigskip
{\bf Abstract}
{\it A general class of finitely ramified fractals 
is that of  P.C.F. self-similar sets. 
An important open problem in analysis on fractals was 
whether there exists a self-similar energy
on every P.C.F. self-similar set.
In this paper, I solve the problem, showing an example of a P.C.F. self-similar set
where there exists no self-similar energy.}



\bigskip

\centerline{\bf 1. Introduction.}

\medskip

\qua An important problem in analysis of fractals is 
the construction of a Laplace operator, or equivalently,
an energy, more precisely, a  self-similar Dirichlet form. 
The construction of a  self-similar Dirichlet form
has been investigated  specially on     finitely ramified fractals. 
Roughly speaking, a fractal is  finitely ramified
 if the intersection of each pair of copies
of the fractal is a finite set. The Sierpinski Gasket, the 
Vicsek Set and the Lindstr\o m Snowflake are finitely 
ramified fractals, while the Sierpinski Carpet is not.

\qua More specially, we consider the 
P.C.F. self similar-sets,  a general  class
of finitely ramified fractals  introduced by Kigami in [3]. 
 A  general theory with many examples can be found in [4]. 
On such a class of fractals,
 the basic tool used to construct a self-similar Dirichlet form
 is a {\it discrete}
 Dirichlet form defined on a special
 finite subset $V^{(0)}$ of the fractal. Such discrete Dirichlet forms 
 have to be
 {\it  eigenforms}, i.e. the eigenvectors of a special
  nonlinear operator $\Lambda_r$ called {\it renormalization operator}, 
  which depends on a set of positive  {\it weights} $r_i$ placed on the cells of the fractal.
In [5], [10]  and  [6] criteria for the existence 
of an eigenform with prescribed weights are discussed. In particular, in [5], 
T. Lindstr\o m proved that there exists
an eigenform on the nested fractals with all weights equal to $1$,
 C. Sabot in [10] proved a rather general criterion, and V. Metz in [6]  
 improved the results in [10].

\qua  In [1],  [7],   [8] and
[9], instead, the problem is considered
whether  on a given fractal   there exists  a {\it G-eigenform}.
By this we mean a form $E$ which is an eigenform
of  the operator $\Lambda_r$ for some set of weights $r$.
In such papers the existence of a G-eigenform was proved on some
classes of P.C.F. self-similar sets.
In fact,  the following open problem is well-known

\centerline{\it  Does  a G-eigenform exist on every 
P.C.F.self-similar set?}

In this paper, I solve such  a problem showing an example
of P.C.F.self-similar set with no G-eigenform.
Here, I consider a very general class of P.C.F.
self-similar sets, as usually in my previous papers on this 
topic (see Section 2 for the details), but also in other papers
(for example it is considered in  [11].
The example is constructed in Section 3. It is a variant
of the $N$-Gaskets, in the sense that every cell $V_i$ only intesects
$V_{i-1}$ and $V_{i+1}$, and has twenty vertices. The ptoof
is based on the evaluation of the effective conductivities
on pairs of close vertices and of far vertices. 
Note that in [9] the existence of a G-eigenform is 
proved on every fractal (of the class considered here)
but only if we consider the fractal generated
by a set of similarities which is not necessarily the given set of similarities
(see [9] for the details).

\bigskip

\centerline{\bf 2. Definitions and Notation.}

\medskip
I will now   define the fractal setting, which is based on that in [9]. 
 This kind of approach was firstly given in 	[2].
We define a fractal  by giving a {\it fractal triple}, i.e., a triple
$\cf:=(V^{(0)},V^{(1)},\Psi)$ where $V^{(0)}=V$ and $V^{(1)}$ 
are finite sets with $\# V^{(0)}\ge 2$,
and $\Psi$ is a finite set of one-to-one maps from $V^{(0)}$ into
$V^{(1)}$ satisfying 

$$V^{(1)}=\bigcup\limits_{\psi\in\Psi}\psi(V^{(0)})\, .$$

  We put $V^{(0)}=\big\{P_1,...,P_N\big\}$, and of course
$N\ge 2$. A set of the form $\psi(V^{(0)})$ with $\psi\in\Psi$
will be called a cell or a $1$-cell.
We require that 

\smallskip
\qua {\it a) For each $j=1,...,N$ there exists a (unique) map $\psi_j\in\Psi$ 
such that $\psi_j(P_j)=P_j$, and
$\Psi=\big\{\psi_1,...,\psi_k\big\}$, with $k\ge N$. 

\smallskip
\qua b) $P_j\notin\psi_i(V^{(0)})$ when $i\ne j$ (in other words, if $\psi_i(P_h)=P_j$ with
$i=1,...,k$, $j,h=1,...,N$, then $i=j=h$).

\smallskip
\qua c) Any two points
in $V^{(1)}$ can be connected by a path  any edge of which
 belongs to a $1$-cell, depending 
of the edge.}
\smallskip

Of course, it immediately follows $V^{(0)}\subseteq V^{(1)}$.
Let $W=]0,+\infty[^k$ and  put
$V_i=\psi_i(V^{(0)})$ for each $i=1,...,k$.   
Let $J(=J(V))=\{\{j_1,j_2\}:j_1,j_2=1,...,N, j_1\ne j_2\}$.
It is well-known that on every fractal triple we can construct
a P.C.F.-self-similar set.

\qua  We denote by $\cd(\cf)$ or simply $\cd$ the set
of the  Dirichlet forms on $V$, invariant with respect to
an additive constant,  i.e., the set of the functionals $E$ 
from $\bre^{V}$ into $\bre$ of the form
$$E(u)=\sum\limits_{\{j_1,j_2\}\in J}
  c_{\{j_1,j_2\}}(E)\big(u(P_{j_1})-u(P_{j_2})\big)^2 $$
with $c_{\{j_1,j_2\}}(E)\ge 0$. I will denote by
$\td(\cf)$ or simply $\td$ the set of the irreducible Dirichlet forms,
 i.e., $E\in\td$ if $E\in
\cd$ and moreover  $E(u)=0$
 if and only if $u$ is constant.
The numbers $c_{\{j_1,j_2\}}(E)$ are called {\it coefficients}
 of $E$. We also say that $c_{\{j_1,j_2\}}(E)$ is the {\it conductivity}
 between $P_{j_1}$ and $P_{j_2}$ (with respect to $E$).
 Next, I recall the notion of effective conductivity. 
Let $E\in\td$, and let $j_1,j_2=1,...,N$, $j_1\ne j_2$. Then we put
$$\cl_{V; j_1, j_2}=\big\{u\in\bre^V: u(P_{j_1})=0,
 u(P_{j_2})=1\big\}.$$
It can be easily proved that the minimum
$\min\{E(u):u\in \cl_{V; j_1, j_2}\}$ exists, is attained at a unique function,
and amounts to $\min\{E(u):u\in \cl_{V; j_2, j_1}\}$.
So, for $E\in\td$ and  $\{j_1,j_2\}\in J$, we define 
$\tc_{\{ j_1,j_2\}}(E)(=\tc_{ j_1,j_2}  (E))$ by
$$\tc_{\{ j_1,j_2\}}(E)=
\min\{E(u):u\in \cl_{V; j_1, j_2}\}=\min\{E(u):u\in \cl_{V; j_2, j_1}\}.$$
The value $\tc_{\{ j_1,j_2\}}(E)$ or short $\tc_{\{ j_1,j_2\}}$,
is called {\it effective conductivity}
 between $P_{j_1}$ and $P_{j_2}$ (with respect to $E$).
 Note that $\tc_{\{ j_1,j_2\}}>0$.
The following remark can be easily verified (see Remark 2.9 in 
[9].

\bigskip
{\bf Remark 2.1}.
If $j_1, j_2=1,...,N$, $j_1\ne j_2$,  and 
$E\in\td$, then
  $$\min \big\{E(u): u\in\rvo: u(P_{j_1})=t_1, u(P_{j_2})=t_2\big\}
= (t_1-t_2)^2\, \tc_{\{j_1,j_2\}}.$$
\bigskip

\qua Recall that   for every  $r\in W:= ]0,+\infty[^k$,  ($r_i:=r(i))$
 the {\it renormalization operator} is defined as follows:
  for every $E\in\td$ and every $u\in\rvo$,  
    
 $$\Lb_{r}(E)(u)=\inf\Big\{S_{1,r}(E)(v), v\in\cl(u)\Big\},$$
 $$S_{1,r}(E)(v):=\sum\limits_{i=1}^k r_i E(v\circ\psi_{i}),\quad 
\cl(u):= \big\{v\in \bre^{V^{(1)}}: v=u
  \ \text{on}\ V^{(0)}\big\}.$$
 
 It is well known that $\Lb_r(E)\in\td$ and that
 the infimum is attained at a unique function $v:=H_{1,E;r}(u)$. 
 When $r\in W$, 
 an element $E$ of $\td$ is said to be an {\it $r$-eigenform} with 
 {\it eigenvalue} $\rho>0$ 
 if $\Lb_{r}(E)=\rho  E$. As this amounts to
 $\Lb_{{r\over \rho}}(E)= E$,  we could also require $\rho=1$.
The problem discussed in the present paper is that
of the existence of a 
 G-eigenform in $\td$, in other words, the existence of $E\in\td$ such that 
$\Lb_{r}(E)=\rho E$ for some $\rho>0$ and $r\in W$. 
In next section, I will describe an example of a fractal triple where
there exists no G-eigenform.  To this aim, it will be useful the following 
standard lemma  (see e.g., Lemma 3.3 in [9].

\bigskip
{\bf Lemma 2.2} {\sl For every $E\in\td$ and $\{j_1,j_2\}\in J$
we have
\smallskip

\centerline{$\tc_{\{j_1,j_2\}}\big(\Lb_{r}(E)\big)=
\min\big\{S_{1,r}(E)(v): v\in  H_{j_1,j_2} \big\},$}
\centerline{$\text{where}\ \  H_{j_1,j_2}=\big\{v\in\rvon: v(P_{j_1})=0,v(P_{j_2})=1\big\}\, .
 $}}

\bigskip
\centerline{\bf 3. The Example.}

\medskip
Let $\bar\cf=(\vo,\von,\Psi)$ be a fractal triple so defined. Let $N=2N'$
 be a positive even number.
Let $\vo=\{P_1,...,P_N\}$, and we fix $N'=10$ so $N=20$.
Let $\Psi=\{\psi_1,...,\psi_{20}\}$. Here, thus, $N=k=20$. 
In the following, the indices of the points and of the maps
will be meant to be mod $20$. For example, $i+9=2$ if $i=13$.
Suppose $V_i \cap V_{i'}=\varnothing$ if $i'\notin \{ i, i-1,i+1\}$, and

$$V_i\cap V_{i+1}=\{\bar Q_i\},\quad \bar Q_i=\cases
\psi_i(P_{i+1})=\psi_{i+1}(P_i)\quad\  \text{\ if\ } i\text{\  odd,}\\
\psi_i(P_{i+9})=\psi_{i+1}(P_{i-8})\ \ \text{\ if\ } i\text{\  even.}\endcases
$$
In this way $\bar Q_{i-1}=
\psi_i(P_{\sigma(i)})$, $\bar Q_i=\psi_i(P_{\sigma(i)+10})$, with
$\sigma(i)=\cases
i-1\text{\ if} \ i\text{\  even,}\\
i-9    \text{\ if\ } i\text{\  odd}. \endcases$
Thus the points $\bar Q_i$ and $\bar Q_{i-1}$ are opposite in $V_i$.
Here we say that $P_h$ and $P_{h+10}$ are opposite in $\vo$ and that
$\psi_i(P_h)$ and $\psi_i(P_{h+10})$ are opposite in $V_i$.
In order to prove Theorem 3.2, we could use arguments based on
effective resistances in series, but, in order to avoid some slightly
technical points, I prefer to give a direct proof. 
We need the following well known lemma.

\bigskip
{\bf Lemma 3.1.}
{\sl For every positive integer $n$  and every $b_i>0$, $i=1,...,n$, 
we have
$$\Big(\sum\limits_{i=1}^n b_i^{-1}\Big)^{-1}=
\min\Big\{\sum\limits_{i=1}^n b_i x_i^2:x\in Y_n\Big\},
\quad Y_n:=\Big\{x\in \bre^n:\sum \limits_{i=1}^nx_i=1\Big\}.
$$}

\smad
Let $f:\bre^n\to\bre$ be defined by $f(x)=\sum\limits_{i=1}^n b_i x_i^2$.
Since $f$ is continuous and $f(x)\llim_{x\to \infty} +\infty$,
$f$ attains a minimum $m$ on the closed set $Y_n$ at some point $\bar x$.
We find $\bar x$ using the Lagrange multiplier rule.
We have $b_i \bar x_i= \lambda$ for some $\lambda\in\bre$ and every $i=1,...,n$.
Thus,
$1=\sum\limits_{i=1}^n \bar x_i=\lambda\sum\limits_{i=1}^n b_i^{-1}, 
$
and $\bar x_i=\lambda b_i^{-1}=b_i^{-1}
\Big(\sum\limits_{j=1}^n b_j^{-1}\Big)^{-1}$. Since $m=f(\bar x)$, 
a simple calculation completes the proof.  \enpr

\bigskip
{\bf Theorem 3.2.}
{\sl On $\bar\cf$ there exists no $G$-eigenform.}

\smad
Suppose by contradiction there exist $E\in\td$ and $r\in W$
such that $\Lb_r(E)=E$. Of course, in view of Lemma 2.2, this implies
$$\tc_{`\{j_1,j_2\} }(E)=\min\big\{S_{1,r}(E)(v): v\in H_{j_1,j_2}\big\} 
\quad \forall\, \{j_1,j_2\}\in J.\eqno (3.1)$$
Now, let 
$\bar r=\max\Big\{\min\{r_{2h+1},r_{2h+2}\}: h=0,...,9\Big\}$.
Thus, 
$$\exists\, \bar h=0,...,9: \forall\, d=1,2: r_{2 \bar h+d} \ge \bar r, \ 
 \eqno (3.2)$$
$$\forall\, h=0,...,9 \ \exists\, d=1,2: r_{2 h+d}\le \bar r. \eqno (3.3)
$$
Next, we evaluate $\tc_{2\bar h+1, 2\bar h+2}(E)$ using (3.1).
Let $v\in H_{2\bar h+1, 2\bar h+2}$.
Then, since by definition $S_{1,r}(E)(v)=\sum\limits_{i=1}^{20} r_iE(v\circ\psi_i)
\ge \sum\limits_{d=1}^{2} r_{2\bar h+d} E(v\circ\psi_{2\bar h+d})$,
in view of (3.2) we have 
$$S_{1,r}(E)(v)\ge 
 \bar r\big(E(v\circ \psi_{2\bar h+1})+E(v\circ \psi_{2\bar h+2})\big).\eqno (3.4)$$
  Let $t=v(\bar Q_{2\bar h+1})$. Let $v_1=v\circ \psi_{2\bar h+1}$,
  $v_2=v\circ \psi_{2\bar h+2}$. Then we have 
  $$v_1(P_{2\bar h+1})=v({P_{2\bar h+1}})=0,\quad
 v_1(P_{2\bar h+2})=v\big(\psi_{2\bar h+1}(P_{2\bar h+2})\big)=v(\bar Q_{2\bar h+1})=t,$$
  $$v_2(P_{2\bar h+2})=v({P_{2\bar h+2}})=1,\quad 
 v_2(P_{2\bar h+1})=v\big(\psi_{2\bar h+2}(P_{2\bar h+1})\big)=v(\bar Q_{2\bar h+1})=t.$$
  By Remark 2.1 and Lemma 3.1 with $n=2$,
  we thus have 
  $$E(v\circ \psi_{2\bar h+1})+E(v\circ \psi_{2\bar h+2})\ge
 \big (t^2+(1-t)^2\big) \tc_{2\bar h+1, 2\bar h+2}(E)
 \ge {1\over 2}\tc_{2\bar h+1, 2\bar h+2}(E). $$
 By (3.4) and (3.1) we have 
 $\tc_{2\bar h+1, 2\bar h+2}(E)\ge \disp{{\bar r\over 2} }\tc_{2\bar h+1, 2\bar h+2}(E)$,
 thus
 $$\bar r\le 2. \eqno (3.5)$$
 Let now $\whh=1,...,20$ be so that
 $$\tc_{\whh,\whh+10}(E)\ge \tc_{l,l+10}(E)\quad \forall\, l=1,...,20. \eqno (3.6)$$
 Note that, in view of Remark 2.1, for every $i=1,...,20$ and every $t_1,t_2\in\bre$ there
 exists  $u_{i,t_1,t_2}\in \rvo$ such that
$$u_{i,t_1,t_2}(P_i)=t_1,\ u_{i,t_1,t_2}(P_{i+10})=t_2,
\ E(u_{i,t_1,t_2})=(t_1-t_2)^2 \, \tc_{i,i+10}(E). \eqno(3.7)$$
Next, define $v\in H_{\whh,\whh+10}$, in terms of $x,x'\in Y_9$.
Given such $x,x'$ let $s(n)=\sum\limits_{i=1}^n x_i$, $s'(n)=1- 
\sum\limits_{i=1}^n x'_i$ for $n=0,...,9$. Let $v\in\rvon$ be so that
$$
v\circ\psi_i=\cases 0\quad\quad\quad\quad\quad\quad\quad\quad\,\  \
\quad\quad \text{if}\ i=\whh\\
1\quad\quad\quad \quad\quad\quad\quad\quad\quad\quad\,\ 
\ \text{if}\ i=\whh+10\\ 
u_{\sigma(i), s(i-\whh-1), s(i-\whh)}\quad\quad\ \ \ \, \ \text{if}\ i=
\whh+1,...,\whh+9,\\ 
u_{\sigma(i), s'(i-\whh-11), s'(i-\whh-10)}\quad\, \text{if}\ i=
\whh+11,...,\whh+19.\\ 
\endcases 
$$
We easily see that the definition of $v$ is correct, i.e., the definition
of $v$ at the points $\bar Q_i$ (the only points lying in different cells)
is independent of the two representations of $\bar Q_i$, that is,
$\bar Q_i=\psi_i(P_{\sigma(i)+10})$ and $\bar Q_i=\psi_{i+1}(P_{\sigma(i+1)})$.
 Moreover, as $P_{\whh}=\psi_{\whh}(P_{\whh})$, and
 $P_{\whh+10}=\psi_{\whh+10}(P_{\whh+10})$.
 we immediately see that $v\in H_{\whh,\whh+10}$. 
 Since $E(v\circ\psi_{\whh})=E(v\circ\psi_{\whh+10})=0$, we have
 $$S_{1,r}(E)(v)=\sum\limits_{i=\whh}^{\whh+19} r_i E(v\circ\psi_i)=
 \Big(\sum\limits_{i=\whh+1}^{\whh+9} r_i E(v\circ\psi_i)+
 \sum\limits_{i=\whh+11}^{\whh+19} r_i E(v\circ\psi_i)\Big)$$
 $$=\sum\limits_{i=\whh+1}^{\whh+9} r_i
 E\big(u_{\sigma(i), s(i-\whh-1), s(i-\whh)}   \big)+
 \sum\limits_{i=\whh+11}^{\whh+19} r_i 
 E\big( u_{\sigma(i), s'(i-\whh-11), s'(i-\whh-10)}\big)$$
 $$
=\sum\limits_{i=\whh+1}^{\whh+9} 
r_i \ \tc_{\sigma (i), \sigma (i+10)} \  x_{i-\whh}^2+
 \sum\limits_{i=\whh+11}^{\whh+19}
 r_i \  \tc_{\sigma (i), \sigma (i+10)}\ (x'_{i-\whh-10})^2 \eqno \text{by}\  (3.7)$$
 $$=\sum\limits_{i=1}^{9}r_{i+\whh}  
 \ \tc_{\sigma (i+\whh), \sigma (i+\whh)+10} \ x_{i }^2
 +\sum\limits_{i=1}^{9}    r_{i+\whh+10} \tc_{\sigma (i+\whh+10),
  \sigma (i+\whh+10)+10} \ (x'_{i })^2$$
  $$\le 
  \tc_{\whh,\whh+10}(E)\Big(
  \sum\limits_{i=1}^{9}r_{i+\whh }  \  x_{i }^2+
  \sum\limits_{i=1}^{9}    r_{i+\whh+10} (x'_{i })^2\Big). \eqno \text{by}\  (3.6)
  $$ 
So far we have taken arbitrary $x, x'\in Y_9$. Now take those
$x,x'\in Y_9$ that minimize the sums in previous formula.
By Lemma 3.1, with this $v$ we have
$$
S_{1,r}(E)(v)\le \tc_{\whh,\whh+10}(E)\Big(
 \big( \sum\limits_{i=1}^{9}(r_{i+\whh})^{-1} \big)^{-1}  +
  \big(\sum\limits_{i=1}^{9}    (r_{i+\whh+10})^{-1} \big)^{-1}\Big).
  \eqno (3.8)
$$
By (3.3), for at least four $i=1,...,9$ we have $r_{i+\whh}\le \bar r$,
so $\big(\sum\limits_{i=1}^{9}(r_{i+\whh})^{-1} \big)^{-1} <{\bar r\over 4}$.
Similarly, $\big(\sum\limits_{i=1}^{9}    (r_{i+\whh+10})^{-1} \big)^{-1} 
<{\bar r\over 4}$. Since $v\in H_{\whh,\whh+10}$, by (3.8) and (3.1)
we have
$$\tc_{\whh,\whh+10} (E) \le S_{1,r}(E)(v)
<  {\bar r\over 2}   \tc_{\whh,\whh+10} (E).$$
Thus, $\bar r>2$, which contradicts (3.5). \enpr

\bigskip

\centerline{\bf References}

\medskip
[1] B.M. Hambly, V. Metz, A. Teplyaev, {\it 
Self-similar energies on post-critically finite self-similar fractals},
 J. London Math Soc. (2) 74, pp. 93-112, 2006

\smallskip
[2] K.Hattori, T. Hattori, H. Watanabe, {\it Gaussian field theories 
on general networks and the spectral dimension}, 
Progr. Theoret. Phys. Suppl. 92, pp. 108-143, 1987

\smallskip
[3] J. Kigami, {\it Harmonic calculus on p.c.f. self-similar sets}, Trans. Amer. Math. Soc. 335, pp.721-755, 1993

\smallskip
[4] J. Kigami, {\it Analysis on fractals}, Cambridge University Press, 2001

\smallskip
[5] T. Lindstr$\ddot o$m,  {\it Brownian motion on nested fractals}, 
Mem. Amer. Math. Soc. 83 No. 420, 1990

\smallskip
[6] V. Metz,  {\it The short-cut test}, J. Funct. Anal. 220, pp. 118-156, 2005


\smallskip
[7] R. Peirone,  {\it  Existence of eigenforms on fractals with three vertices}, 
Proc. Royal Soc. Edinburgh Sect. A 137, 2007

\smallskip
[8] R. Peirone, \ {\it  Existence of 
eigenforms on nicely separated fractals}, in \ {\it Analysis of graphs and its applications}, Amer. Math. Soc., Providence, pp. 231-241, 2008

\smallskip
[9] R. Peirone, \ {\it  Existence of self-similar energies on
 finitely ramified fractals}, Journal d'Analyse Mathématique 
 Volume 123 Issue 1, pp. 35-94, 2014

\smallskip
[10] C. Sabot,  {\it  Existence and uniqueness of diffusions on 
finitely ramified self-similar fractals}, Ann. Sci. \' Ecole Norm. Sup. (4) 30, pp. 605-673, 1997

\smallskip
[11] R.S. Strichartz, {\it  Differential equations on fractals: a tutorial}, 
Princeton University Press, 2006

\end{

\end{document}